\documentclass[draft]{article}

\usepackage[latin2]{inputenc}
\usepackage{amsmath, amsthm}

\newcommand{\PP}{{\mathcal P}}
\newcommand{\PPn}{{\mathcal P}_n}
\newcommand{\RR}{{\bf R}}
\newcommand{\LL}{{\mathcal L}}
\newcommand{\CC}{{\bf C}}
\newcommand{\KK}{{\bf K}}

\newcommand{\NN}{{\bf N}}

\newcommand{\intt}{{\rm int\,}}
\newcommand{\con}{{\rm con\,}}

\newcommand{\vol}{{\rm vol\,}}

\def\eqref#1{(\ref{#1})}

\newtheorem{theorem}{Theorem}

\newtheorem{corollary}{Corollary}

\newtheorem{proposition}{Proposition}
\newtheorem{definition}{Definition}
\newtheorem{conjecture}{Conjecture}

\newcounter{oldlemma}
\def\theoldlemma{\Alph{oldlemma}}
\newenvironment{oldlemma}[1][1]{
  \em
  \vskip 1ex
  \refstepcounter{oldlemma}
  \noindent{\bf Lemma\ \theoldlemma\ifx#1 1\relax\else \ \rm(\bf#1)\fi{\bf.}}
}{\vskip 1ex}

\newcounter{oldresult}
\def\theoldresult{\Alph{oldresult}}
\newenvironment{oldresult}[1][1]{
  \em
  \vskip 1ex
  \refstepcounter{oldresult}
  \noindent{\bf Theorem\ \theoldresult\ifx#1 1\relax\else \ \rm(\bf#1)\fi{\bf.}}
}{\vskip 1ex}

\newcounter{oldprop}
\def\theoldprop{\Alph{oldprop}}

\begin{document}

\title{INEQUALITIES FOR MULTIVARIATE POLYNOMIALS}
\author{Szil\'ard Gy. R\'ev\'esz}

\date{}

\maketitle

\begin{abstract}
We summarize researches -- in great deal jointly with my host Y.~Sarantopoulos and his PhD
students V.~Anagnostopoulos and A.~Pappas -- started by a Marie Curie fellowship in 2001 and is
still continuing.

The project was to study multivariate polynomial inequalities. In
the course of work we realized the role of the ``generalized
Minkowski functional'', to which we devoted a throughout survey.
Resulting from this, infinite dimensional extensions of
Chebyshev's extremal problems were tackled successfully.
Investigating Bernstein-Markov constants for homogeneous
polynomials of real normed spaces led us to the application of
potential theory. Also we found at first unexpected connections of
polarization constants of $\RR^2$ and $\CC^2$ to Chebyshev
constants of $S^2$ and $S^3$, respectively. In the study of
polarization constants, a further application of potential theory
occurred. This led us to realize that the theory of rendezvous
numbers can be much better explained by potential theory, too. Our
methods for obtaining Bernstein type pointwise gradient estimates
for polynomials were compared in a recent case study to the yields
of pluripotential theoretic methods. The findings were that the
two rather different methods give exactly the same results, but
the two currently standing conjectures mutually exclude each
other.
\end{abstract}

\section{Polynomials In Higher Dimensions}

In the whole paper $\KK$ stands for either $\RR$ or $\CC$. If
$X$ is a normed space, $X^*=\LL(X,\KK)$ is the usual dual space,
and $S:=S_X$, $S^*:=S_{X^*}$, $B:=B_X$ and $B^*:=B_{X^*}$ are the
unit spheres and (closed) unit balls of $X$ and $X^*$,
respectively. Moreover, $\PP=\PP(X)$ and $\PP_n=\PP_n(X)$ will
denote the space of continuous (i.e., bounded) polynomials of free
degree and of degree $\leq n$, resp., from $X$ to $\KK$.

There are several ways to introduce continuous polynomials over $X$, one being the linear
algebraic way of writing
\begin{equation}\label{Pnassum}
\PP_n:=\PP^*_0+\PP^*_1+\dots+\PP^*_n ~, \quad {\rm and} \quad \PP:=\bigcup_{n=0}^{\infty}\PPn
\end{equation}
with $\PP^*_k:=\PP(^kX;\KK)$ denoting the space of homogeneous (continuous) polynomials of
degree (exactly) $k\in\NN$. That is, one considers bounded $k$-linear forms
$ L\in \LL(X^k\to\KK) $
together with their ``diagonal functions''
\begin{equation}\label{Lhat}
\widehat{L}:X\to\KK, \quad \widehat{L}(x):=L(x,x,\dots,x),
\end{equation}
and defines $\PP^*_k$ as the set of all $\widehat{L}$ for $L$ running
$\LL(^kX):=\LL(X^k\to\KK)$. In fact, it is sufficient to identify equivalent linear forms
(having identical diagonal functions) by selecting the unique symmetric one among them: i.e.,
to let $L$ run over $\LL^s(^kX)$ denoting
{\em symmetric} $k$-linear forms. Building up the notion of polynomials so is equivalent to 
\begin{alignat}{1}\label{PnX}
\PPn:=\{p:& ~ X\to\KK~:~\|p\|<\infty, p|_{Y+y}\in\PPn(\KK)  \\
  &\qquad \quad{\rm for}~{\rm all}~~ Y\le X,~\dim Y =1,~ y\in X \},\nonumber
\end{alignat}
or to the definition arising from combining \eqref{Pnassum} and
\begin{alignat}{1}\label{PkXstar}
\PP^*_k:= \{p : &~X\to\KK~:~\|p\|<\infty, ~p|_{Y} \in\PPn^*(\KK^2) \\
 & \qquad \quad {\rm for}~{\rm  all} ~~ Y\le X,~\dim Y =2 \}.\nonumber
\end{alignat}
Here and throughout the paper for any set $K\subset X$ and function $f:X\to\KK$ we denote, as
usual,
$$
\|f\|_K:=\sup_K |f| \qquad {\rm and}\qquad \|f\|:=\|f\|_B~.
$$
For an introduction to polynomials over normed spaces see, e.g., \cite[Chapter 1]{Din}. In
particular, it is well-known that
\begin{equation}\label{LhatLLhat}
\|\widehat{L}\|\le\|L\|\le C(n,X) \|\widehat{L}\|\quad {{\rm for}~{\rm all}}~~ L\in\LL^s(^n X),
\end{equation}
and that $C(n,X)\le n^n/n!$, \cite{Din}, while $C(n,X)=1$ if $X$ is
a Hilbert space (Banach's Theorem). Similarly to \eqref{LhatLLhat},
one can consider special homogeneous polynomials which can be
written as products of linear forms, that is
$L(x_1,\dots,x_n)=\prod_{j=1}^n f_j(x_j)$ with $f_j\in X^*$. Then
$\|L\|=\prod_{j=1}^n\|f_j\|$, i.e., the product of the norms, and
one compares to the norm of the corresponding homogeneous
polynomial, i.e., to $\|\widehat{L}\|=\|\prod_{j=1}^n f_j\|$. Note
that here $L$ is far from being symmetric, and this yields to an
essentially different question, with the similarly defined
polarization constants -- the \emph{linear polarization constants}
-- now ranging up to $n^n$, see e.g.~\cite{BST}. Polarization
problems are typical, genuinely multivariate inequalities, as in
dimension 1 they degenerate.

\section{Linear Polarization Constants}

\begin{definition}[\bf Ben\'\i tez, Sarantopoulos, Tonge
\cite{BST}\rm]\label{nthpolarconst} The $n^{\rm th}$ {\rm (linear)
polarization constant} of a normed space $X$ is
\begin{alignat}{1}\label{polconstdef}
 c_n(X): &= \inf \{M : \prod\limits_{j=1}^n \|f_{j}\| \leq M \|\prod\limits_{j=1}^n
f_{j}\|\; (\forall f_{j} \in X^{\ast}) \}\nonumber \\
& = 1/ \inf\limits_{f_{1}, \ldots ,f_{n} \in S_{X^{\ast}}} \sup\limits_{\|x\|=1}|\prod_{j=1}^n
f_{j}(x)|\,.
\end{alignat}
\end{definition}

\noindent Obviously $c_{n}(X)$ is a nondecreasing sequence. Its growth, as $n \rightarrow
\infty$, is closely related to the structure of the space.

\begin{definition}[\bf R\'ev\'esz, Sarantopoulos \cite{KRS}\rm]
The {\rm linear polarization constant} of a normed space $X$ is
\begin{equation}
c(X):=\lim_{n \rightarrow \infty} c_{n} (X)^{ \frac{1}{n}}\;\;.
\end{equation}
\end{definition}
\noindent One should have put only $c(X):=\limsup_{n \rightarrow \infty} c_{n} (X)^{
\frac{1}{n}}$ as a definition. However, we proved that the limit does exist, see
\cite[Proposition 4]{KRS}. Note that $c(X)$ can be infinity as well. More specifically, from
\cite[Theorem 12]{KRS} we have

\begin{proposition}[\bf R\'ev\'esz, Sarantopoulos \cite{KRS}\rm]\label{infdim}
Let $X$ be a normed space. Then $c(X)=\infty$ iff $\dim (X)=\infty$.
\end{proposition}

In the special case where $X$ is a Hilbert space, it is easy to see that (writing $Y\leq X$ for
$Y$ being a subspace of $X$)
\begin{equation}\label{finite}
c_{n}(X)=\sup \left\{c_{n} (Y): Y\leq X 
\,, \dim Y\le n \right\}\,.
\end{equation}

The Banach-Mazur distance $d(X, Y)$ between two isomorphic Banach spaces $X$ and $Y$ can be
used in comparing the $n^{\rm th}$ polarization constants of these spaces. Recall that
\[
d(X,Y):= \inf \{ \|T\| \cdot \|T^{-1}\|:\, T : X \leftrightarrow Y\: {\mbox {isomorphism}}
\}\,.
\]

\begin{proposition}[\bf Ben\'\i tez, Sarantopoulos, Tonge
\cite{BST}\rm]\label{BeSaTo} If the normed spaces $X$ and $Y$ are isomorphic, then
\begin{equation}\label{BanachMazur}
c_n(X) \leq \ d^n(X, Y) c_n(Y)\,.
\end{equation}
\end{proposition}
\noindent It seems very likely that $c_n(X) \geq c_n(\ell_2^n)\,\,(\forall n \le \dim X)$, but
this is not known. However, we found
\begin{proposition}[\bf R\'ev\'esz, Sarantopoulos \cite{KRS}\rm]\label{Kadets}\hskip0.5em
If $X$ is an infinite dimensional normed space, then
\begin{equation}\label{polarconstant}
c_n(X) \geq c_n(\ell_2^n)\,,\forall {n \in \NN} \,.
\end{equation}
\end{proposition}
As is well-known, for any $n$-dimensional Banach space $X$
\begin{equation}\label{John} d(X, \ell_{2}^n)\leq \sqrt{n}\,.
\end{equation}
Thus to determine $c_{n}({\bf R}^{n})$ is interesting not only in the context of Hilbert space
theory. For example, a combination of \eqref{finite}, \eqref{BanachMazur}, \eqref{John} and
\eqref{polarconstant} yields the following result.
\begin{theorem}[\bf R\'ev\'esz, Sarantopoulos \cite{KRS}\rm]\hskip-0.1em
Let $X$ be an infinite dimensional normed space and let $H$ be the space $\ell_{2}$ over $
\KK$. Consider $\KK^{n}$, i.e., the space $\ell_2^n$ over $\KK$. For all $n \in \NN $ we have
\begin{equation*}
c_{n}(H) = c_{n}( {\KK}^{n}) \leq c_{n}(X) \leq n^{\frac n2} c_{n} ( {\KK}^{n}) = n^{\frac n2}
c_{n}(H)\,.
\end{equation*}
\end{theorem}
Note that determination of the linear polarization constant is closely connected to another
famous problem, the Tarski plank problem, see, e.g., \cite{BALL2, KRS}.

Let us focus here on the problem of estimating the linear polarization constants of Hilbert spaces. Although this is a classical topic, 
there was a flourishing activity on this field even in the last ten years, 
and even recently, see, e.g., \cite{BALL2, Mate, MS}.

\vskip1ex\noindent Ideally, one should look for the exact values of $c_n(\ell_2^d)=c_n(\KK^d)$,
for any $d, n \in \NN$, which, in view of \eqref{finite}, reduces to $d\le n \in \NN$. In fact,
this question is posed in \cite{LMS}, attributed to the referee of the paper. In this direction
a remarkable success is Arias-de-Reyna's result.

\begin{oldresult}[\bf Arias-de-Reyna \cite{ARIAS}\rm]\label{Ar}
$ c_{n}( {\CC}^{n})=n^{n/2}\;\;. $
\end{oldresult}

An even more precise description was obtained recently by K. Ball \cite{BALL2}. Estimating
$c_{n}( \bf{R}^{n})$ seems to be a harder task. In particular, the proof of the stronger result
due to K. Ball relies heavily on complex function theoretic tools which are not valid in the
case of $ \bf{R}^{n}$. Observe furthermore that Arias-de-Reyna's entirely different technique,
based on per\-manents, multilinear algebra and probability theory (particularly Gaus\-sian
random variables), strongly depends on the complex structure of $\bf{C}^{n}$. As Arias-de-Reyna
has mentioned in \cite{ARIAS}, his Theorem \ref{Ar} leads to an upper estimate $c_n(\RR^n)\leq
2^n n^{n/2}$ even for the real case. This has been improved in \cite{LMS} and \cite{GV}, until
a more refined approach was worked out in \cite{MST}, using the natural complexification of a
real Hilbert space. This has led to the currently best

\begin{theorem}[\bf R\'ev\'esz, Sarantopoulos \cite{KRS}\rm]\label{Resa}
$$
n^{ \frac{n}{2}} \leq c_{n} ( \RR^{n}) \leq 2^{ \frac{n}{2}-1} n^{ \frac{n}{2}}\;\;.
$$
\end{theorem}

Let us mention here the following conjecture, appearing already in \cite{BST} and \cite{ARIAS}
and formulated also in \cite{KRS}. \vskip 0.10in
\begin{conjecture}
\begin{equation}\label{Sarconj}
c_{n}( \RR^{n} )=n^{ \frac{n}{2}}\;\;.
\end{equation}
\end{conjecture}
We proved the conjecture for $n=1,2,3,4,5$ in \cite{PR}. Moreover, in \cite{PR} we discussed a
direct, real approach in detail. This approach seems to be interesting (even though the
resulting exponential factor falls, unfortunately, only between $\sqrt 2$ and $2$), as it is
independent of Theorem \ref{Ar}.

In all, for a Hilbert space $H$ of infinite dimension we only know that $c(H)=\infty$ and $\log
c_{n} (H) \sim \frac{1}{2}n \log n$. On the other hand, if $\dim H=d$ is fixed then by Theorem
\ref{infdim} $c(H)$ must be finite.

Trying to determine $c_n(\KK^d)$ for arbitrary $d\le n \in \NN$, by natural extrapolation one
might have thought that $c_n(\KK^d)=d^{n/2}$. This was disproved first in \cite{VASILIS}.
\begin{definition}
The $n^{th}$ {\rm (metric) Chebyshev constant} of a subset $F \subseteq X$ in a metric space
$(X, \rho)$ is
\begin{equation}
M_{n} (F):=\inf_{y_{1}, \ldots ,y_{n} \in F} \sup_{y \in F} \rho (y,y_{1}) \cdot \ldots \cdot
\rho(y,y_{n})\;.
\end{equation}
\end{definition}
In particular, in a normed space $X$ $\rho(x,y)=\|x-y\|$ and
\begin{equation}
M_{n} (F):=\inf_{y_{1}, \ldots ,y_{n} \in F} \sup_{y \in F} \|y-y_{1}\| \cdot \ldots \cdot
\|y-y_{n}\| \;.
\end{equation}

\begin{proposition}[\bf Anagnostopoulos, R\'ev\'esz \cite{VASILIS}\rm]\label{Revesz2}
{\sl For \newline the \emph{real} space $\ell_{2}^2(\RR)$ we have
\[
 c_n(\ell_{2}^2(\RR))= \frac{2^n}{ M_n(S^1)}= 2^{n- 1},\:\:\:\mbox{and so}\quad c(\ell_{2}^2(\RR))= 2\,.
\]
Furthermore, for the {\it complex} space $\ell_{2}^2(\CC)$ we have
\[
c_n(\ell_{2}^2(\CC))= \frac{2^n}{ M_n(S^2)},\:\:\:\mbox{and also}\quad
c(\ell_{2}^2(\CC))= \sqrt{e}.
\]}
\end{proposition}

\section{Potential Theory Emerges}

We have seen above that, e.g., $c(\ell_{2}^2(\CC))= \sqrt{e}$, disproving our initial guess of
$c_n(\KK^d)=d^{n/2}$ (which then would have implied $c(\KK^d)=\sqrt{d}$). On the other hand, a
surprising connection with the (metric) Chebyshev constants emerged from our study.

It turned out that for higher dimensions the connection breaks. However, in a more general
setting we still plan to describe linear polarization constants by means of some more general
Chebyshev constants. This is one of the intriguing questions we are occupied with recently.

But how the notion of Chebyshev constants, belonging classically to potential theory, can have
a role here? Quite naturally. Following a potential-theory inspired approach, we could even
describe polarization constants of all finite dimensional Hilbert spaces \cite{PR}. To that,
let us start with a notation:
\begin{equation}\label{Ldef}
L(d,\KK):=\int_{S_{ \KK}^{d}} \log |\langle x, s \rangle | d \sigma (x)\;\; ( <0 )\,,
\end{equation}
where $d \sigma(x)$ is the normalized surface Lebesgue measure of $S_{{\KK }^{d}}$ and $s\in
S_{{\KK }^d}$ is arbitrary. Calculation of the explicit values of the constants $L(d,\KK)$ are
standard.

\begin{oldresult}[\bf Garc\'\i a-V\'azquez, Villa \cite{GV}\rm]\label{spaniards}
We have the eq\-uality $c( \RR^{d})=e^{-L(d,\, {\RR })}$, with the constants $L(d, {\RR })$ as
in \eqref{Ldef}.
\end{oldresult}

\begin{theorem}[\bf Pappas, R\'ev\'esz \cite{PR}\rm]\label{mainresult}For the complex case
we have $c( \CC^{d})=e^{-L(d,\, \CC)}$, with the constants $L(d, \CC)$ defined by \eqref{Ldef}.
\end{theorem}

Our proof in \cite{PR} is a unified, potential theory flavored approach. But the real case were
already obtained a few years earlier by Garc\'\i a-V\'azquez and Villa, with no use of
potential theory at all. Instead, they applied a nice theorem of O. Gross \cite{Gross} on the
existence of \emph{rendezvous numbers}.

Gross' Theorem states that in any compact, connected metric space
$(X,\rho)$ there exists a unique rendezvous number $r=r(X)$, such
that for arbitrary choice of any set of points $x_1,\dots,x_n\in
X$, there always exists some point $x\in X$, such that its average
distance from the $x_j$ is equal to $r$. This is a beautiful
result, which somehow stood in itself for half a century: some
authors even named $r(X)$ ``the magical number'' of $X$.
Nevertheless, several dozens of extensions, applications and
investigations in various contexts were published about rendezvous
numbers.

But how come, that so different approaches can give the same results for polarization
constants? In \cite{FRpotenc} we could satisfactorily describe the theory of rendezvous numbers
by general (linear) potential theory as laid down by, e.g., \cite{Fug}. So it turns out that
the approach of Garc\'\i a-V\'azquez and Villa was not that much different, after all: also
below the surface of Gross' Theorem potential theory lies behind.

\section{More Polynomial Extremal Problems}

The further extremal problems come from natural extensions of classical, univariate
approximation theory questions. Although we do not have space here to describe their widespread
use and various applications, they indeed are rather important questions in particular in
infinite dimensional holomorphy and in multivariate approximation.

One of the main groups of problems may be called Chebyshev type problems of polynomial growth.
A typical, key example is determination of the quantity
\begin{equation}\label{Chebyfindef}
C_n(K, x ):= \sup\big\{p( x ):p\in \PPn,\ {\|p\|}_K\le 1\big\},
\end{equation}
for arbitrary fixed $x \in X$ and $K\subset X$ a fixed convex body.

In dimension 1 the only convex body is the interval, and for the unit interval $I:=[-1,1]$ the
answer is given for all $x\notin I$ by $T_n(|x|)$, with $T_n$ the Chebyshev polynomials
\begin{equation}\label{Tndef}
T_n(x):=
{\frac{1}{2}}\left\{{(x+\sqrt{x^2-1})}^n +{(x-\sqrt{x^2-1})}^n\right\}
\nonumber\,.
\end{equation}

Another group is derivative estimates. A fundamental property of polynomials is that having
control over the size (of values) of a polynomial (say, on a convex body $K$) automatically
provides some finite bound even on their derivatives. One can seek the best bounds at some
given point $x$, (usually within $K$), or some global, uniform bound on the whole of $K$, say.
Also, we can group these problems according to possible further restrictions on the type of
polynomials, or the type of derivatives we consider: e.g., we can restrict to homogeneous
polynomials, or we can consider tangential derivatives etc. A further direction of research is
dealing with higher order derivatives.

First let us consider here the pointwise or Bernstein problem of maximizing $\|\nabla p (x)\|$
for a given $x\in {\rm int} K$ and among all polynomials with $\|p\|_K \le 1$. In analogy
to the best one dimensional estimates, we can normalize the extremal quantity by the
``Bernstein-Szeg\H o factor'' and thus define the ($n^{\rm th}$) Bernstein constant as
\begin{equation}\label{bernsteinfactor}
B_n(K,x):=\frac 1n \sup\limits_{\mbox{\scriptsize $\deg p \le n,$}\atop\mbox{\scriptsize
$|p(x)|<||p||_K$}}\frac {\|\nabla p(x)\|}{\sqrt{||p||_K^2 - p^{2}({x})}}\,.\nonumber
\end{equation}

Next, let us consider the homogeneous Markov factor for arbitrary (say, centrally symmetric)
convex bodies $K$. As then $K$ generates a corresponding norm $\|\cdot\|_{(K)}$, we can
equivalently seek for the estimation of derivatives of homogeneous polynomials in arbitrary
given norms. So put
\begin{equation}\label{HarrisMarkov}
M_n^{(k)}(K):= \sup\limits_{\|x\|_{(K)}\leq 1,\,\,p \in \PPn^*} \|\widehat{D}^k p(x)\| \,,
\end{equation}
where for given $x$ $\|\widehat{D}^k p(x)\|$ stands for the norm of the diagonal $k$-form --
the $k$-homogeneous polynomial -- of the $k$-linear mapping $D^k p(x)$ of $X^k$, with $D^k$
being $k$-fold differentiation. Taking $M_n^{(k)}:=\sup_K M_n^{(k)}(K)$, already Harris
\cite{Harris} has shown that $M_n^{(k)}=c_m^{(k)}$, where
\begin{alignat}{1}\label{harriseq}
  c_m^{(k)}:=\max\{|p^{(k)}(0)|\colon & p\in\PP, \\
  & |p(t)|\le  (1+|t|)^m(t\in\RR)\} .\nonumber
\end{alignat}

\section{Generalized Minkowski Functional}

If $x$ lies in $\RR$ and we want to quantify its position with respect to the unit interval
$I:=[-1,1]$, say, then it suffices to take $|x|$. In several dimensions a more sophisticated
quantitative notion is needed. This was found essentially by Rivlin and Shapiro \cite{RiSh};
but it turned out, that the quantity is actually (equivalent to) an older geometric notion,
going back to Minkowski \cite{Mi} and contemporaries. We extended the notion even to
topological vector spaces, and gave a throughout description of its properties, many equivalent
definitions and some of the related interesting problems, see \cite{RS}. Here we recall just
one definition. By convexity, $K$ is the intersection of its ``supporting halfspaces'' $X(K,
v^*)$, and grouping opposite halfspaces to form layers $L(K, v^*):=X(K, v^*)\cap X(K, -v^*)$ we
get
\begin{equation}\label{(2.9)}
K=\bigcap_{ v^*\in S^*}X(K, v^*)=\bigcap_{ v^*\in S^*} L(K, v^*).
\end{equation}
Any layer can be homothetically dilated with quotient $\lambda \ge 0$ at any of its symmetry
centers lying on its central symmetry hyperplane to get $L^{\lambda}$, and we can even define
\begin{equation}\label{(2.12)}
K^\lambda:= \bigcap_{ v^*\in S^*} L^\lambda(K, v^*).
\end{equation}
Using the convex, closed, bounded, increasing and (as easily seen, cf.~\cite[Proposition 3.3
]{RS}) even absorbing set system ${\{K^\lambda\}}_{\lambda\ge 0}$, the {\it generalized
Minkowski functional\/} is
\begin{equation}\label{alphadef}
\alpha(K, x):= \inf\{\lambda\ge 0\colon x\in K^\lambda\}.
\end{equation}

One advantage of this formulation is that it defines $\alpha$ in a unified way both for $x\in
K$ and for $x\in X\setminus K$. For the Chebyshev problem we need it only for exterior points:
for the Bernstein problem we are concerned with interior points only. The same generalized
Minkowski functional occurs naturally in different extremal problems.

\section{Inequalities For Polynomials}

In what follows, denote $h(K,u^*):=\sup_K u^*$ the \emph{support
functional} and \linebreak $w(K,u^*):=h(K, u^*)+h(K,- u^*)$ the
\emph{width} of $K$ in direction of $u^*\in S^*$. Also, take
$w(K):=\min_{S^*} w(K,\cdot)$ to be the \emph{minimal width} of $K$.
The Chebyshev problem \eqref{Chebyfindef} has the following answer.
\begin{theorem}[\bf Révész, Sarantopoulos \cite{RS}\rm]\label{th:RS}
For an arbit\-ra\-ry convex body $K\subset X$ and any point $ x\notin K$ we have
\begin{equation}\label{CnTn}
C_n(K, x )=T_n(\alpha(K, x )).
\end{equation}
Moreover, $C_n(K, x )$ is actually a maximum, attained by
\begin{equation}\label{Chebyextremal}
P( x ):= T_n \left( \frac{2 \langle v^*, x \rangle -h(K, v^*)+h(K,- v^*)}{w(K, v^*)} \right)
\,.\nonumber
\end{equation}
\end{theorem}
The result was obtained for finite dimensional spaces and strictly convex bodies by Rivlin and
Shapiro \cite{RiSh}. Removing strict convexity is not too difficult, and can be done several
ways, but to extend to infinite dimensions turned to be more delicate. The reason is that not
only compactness is lost, but also the existence of parallel supporting hyperplanes of a
certain configuration. Even though counterexamples show that we no longer have that
configuration, a new proof goes through -- by means of more careful analysis of $\alpha(K,x)$.
For details see \cite{RS, YRS, complutense}.

Even if this seems to settle the question satisfactorily, for many
further applications, in particular for Bernstein and Markov
problems, an extended answer for \emph{complex} points $z\in Z$,
where $Z:=X+iX$ is a \hbox{complexification,} would be very
interesting even for finite dimensional $X$.

Note that even for dimension 1 this is far from being obvious. Indeed,  for complex points
close to $(-1,1)$, Chebyshev polynomials may have values arbitrarily close to 0, thus we can
not expect just one polynomial to be extremal for all $z\in \CC$. On the other hand, results in
this direction can be applied in the Bernstein and Markov problems, see, e.g., \cite{Bojanov}.
For the multivariate complex case we only have some vague ideas at the moment.

\begin{theorem}[\bf R\'ev\'esz, Sarantopoulos \cite{HRS}\rm]\label{th:homogeneous}
With some absolute constant $c_1>0$ we have
\begin{equation}
c_1 m\log m\ \le c_m^{(1)} \le 3m\log m.
\end{equation}
\end{theorem}
The upper estimate was already obtained by Harris \cite{Harris} with some less precise constant
-- the more difficult lower estimate was new. By iterating the result for $M_n^{(k)}$, one even
gets some upper estimate for all $M_n^{(k)}=c_n^{(k)}$. We also computed some better values for
higher derivatives, see \cite{HRS}.

Finally let us consider the Bernstein problem. Here the first, and one of the still very
successful methods -- the method of inscribed ellipses -- were introduced by Sarantopoulos in
1991 \cite{Sar}. The key of all of the method is the \emph{Inscribed Ellipse Lemma:}
\begin{oldlemma}[\bf Sarantopoulos \cite{Sar}\rm]\label{inellipse} Let $K$ be any subset in a vector space $X$. Suppose
that $x \in K$ and the ellipse
\begin{equation}\label{ellipse}
{\bf r}(t) = \cos{t}~ a + b \sin{t}~ y + x-a \qquad (t \in [-\pi,\pi))\,
\end{equation}
lies inside $K$. Then we have for any polynomial $p$ of degree at most $n$ the Bernstein type
inequality
\begin{equation}\label{Bernsteinellipse}
| \langle {D} p(x), y \rangle| \leq \frac{n}{b} \sqrt{ ||p||_K^2 - p^{2}(x)}.
\end{equation}
\end{oldlemma}
The method was applied even to nonsymmetric convex bodies, but in this case our result is still
not final.
\begin{oldresult}[Kro\'o, R\'ev\'esz, Sarantopoulos \cite{KR}, \cite{RS}\rm] \label{th:Bernstein}
Let $K$ be an arbitrary convex body, $x \in {\rm int}K $ and $\|{ y}\| = 1$, where $X$ can be
an arbitrary normed space. Then we have
\begin{equation}\label{krry}
| \langle {D}p(x),  y \rangle | \leq \frac {2 n \sqrt{ ||p||_K^2 - p^{2}({ x}) }}
{\tau(K,y) \sqrt{1 - \alpha(K,x)} } ,
\end{equation}
for any polynomial $p$ of degree at most $n$. Here $\tau(K,y):= \sup \{\lambda\colon \exists
x\in K \,\,{\rm such}\,\, {\rm that}\,\, x+\lambda y \in K\}$ stands for the ``maximal'' chord
in direction $y$.
\end{oldresult}
Since for $x\in K$ $\alpha(K,x)\le 1$, and for $\alpha \in (0,1)$ $\frac1{1-\alpha} \leq \frac
{1+\alpha}{1-\alpha^2}$, we also have
\begin{equation}\label{krrgrad}
B_n(K,x) \leq \frac {2 \sqrt 2 n} {w(K) \sqrt{ 1- \alpha^2(K,x)}}\,.
\end{equation}
Note that apart from the $\sqrt{1+\alpha(K,x)} \leq \sqrt{2}$ factor, the estimate gives the
sharp result even in dimension 1. Hence it was natural to conjecture
\begin{conjecture}[\bf R\'ev\'esz, Sarantopoulos \cite{RS}\rm]\label{alphaconj}
\begin{equation}\label{alphaconjvalue}
B_n(K,x)=\frac{2n}{w(K)\sqrt{1- \alpha^2(K,x)}}~.
\end{equation}
\end{conjecture}
Note that \eqref{krry} was only a -- delicate, but not exact -- estimate. Recently we showed
\cite{MR}, that not even in the case of the standard simplex can the method reach Conjecture
\ref{alphaconj}.

\section{Potential Theory Once More}

We have already seen how potential theory plays a role in the linear polarization constant
problem. It is also well-known, see, e.g., \cite{ST}, that the theory of weighted approximation
and univariate orthogonal polynomials with respect to weights can be analyzed via (weighted)
potential theory of the complex plane. Because the original problem is translated to determine
$c_m^{(k)}$ in \eqref{harriseq}, it is not so much surprising, that with a slight extension of
the theory even our multivariate homogeneous polynomial Markov problem could be treated. That
was our approach in \cite{HRS}.

\vskip1ex\noindent However, there is a genuinely multivariate potential theoretic approach to
multivariate polynomial inequalities: plu\-ripotential theory. Just a few years after
Sarantopoulos, M.~Baran \cite{BARAN2} obtained the same results as \cite{Sar}, with the method
extending to other cases as well.

This theory is well described in, e.g., \cite{BLM, Kli}, so here we summarize only very
briefly. The starting point is the Zaharjuta--Siciak extremal function, which is defined, e.g.,
in $Y:=\CC^d$ with respect to a compact set $E\subset Y$ (or $E\subset X:=\RR^d$, say), as
follows: $V_E$ vanishes on $E$, while outside $E$ we have the definition
\begin{equation}\label{ZSfunction}
V_E(z):=\sup \bigl\{\frac{\log|p(z)|}{\deg p}\,:~ 0\neq p \in \PP(Y),~ ||p||_E\le 1 ~ \bigr\}
\end{equation}
For $E\subset X$ one can easily restrict even to $p\in\PP(X)$. Note that $\log|p(z)|$ is a {\em
plurisubharmonic} function (PSH, for short). In case $E$ is some nice set -- e.g., if it is a
convex body -- then already $V_E$ is continuous. However, even in the general case the upper
semicontinuous regularization $V_E^*$ is at least upper semicontinuous, hence locally bounded
for non-pluri\-po\-lar $E$, which we now assume.

The growth of $(1/\deg p)\log|p(z)|$ is at most $\log_{+}|z|+O(1)$. So it is reasonable to
consider the Lelong class:
\begin{equation}\label{Lelong}
\LL (E) := \{ u\in PSH ~:~u|_E \le 0,~ u(z)\le \log_{+}|z|+O(1)\}\nonumber
\end{equation}
and to define
\begin{equation}\label{Lelongextremal}
U_E(z):=\sup \{ u(z)~:~ u\in\LL (E)\}~.
\end{equation}
This function is named the pluricomplex Green function. The Zaharjuta--Siciak Theorem says that
\eqref{Lelongextremal} and \eqref{ZSfunction} are equal, at least as long as $E\subset \CC^d$
is compact.

The extension of the Laplace- and Poisson equations is the so-called complex Monge--Amp\`ere
equation:
\begin{equation}\label{MAeq}
(\partial\overline{\partial}u)^d := d! 4^d \mathrm{det} \left[ \frac{\partial^2 u}{\partial z_j
\overline{\partial} z_k}(z) \right] dV(z),
\end{equation}
where $dV(z)$ is just the usual volume element in $\CC^d$. Due to the work of Bedford and
Taylor, the operator extends, in the appropriate sense, even to the whole set of locally
bounded PSH functions (which includes $V_E^{*}$ for non-pluri\-po\-lar $E$). Therefore, it
makes sense to consider $(\partial\overline{\partial} V_E^{*})^d$, which is then a compactly
supported measure $\lambda_E$ and is called the {\em complex equilibrium measure} of the set
$E$. The support of $\lambda_E$ lies in the polynomial convex hull $\widehat{E}$ of $E$: in
case $E$ is convex, $\widehat{E}=E$ and $V_E^{*}=V_E$: moreover, it is also shown that
$\lambda|_E(\CC^d)=\lambda|_E(\widehat{E})=(2\pi)^d$. Observe that $V_E(z):=\sup_{n\in\NN}
\frac 1n \log|C_n(E,z)|$  which gives rise for any convex body $K$ and $x\in X\setminus K$ to
the formula $V_K(x)=\log\left(\frac 12 \alpha(K,x)+\sqrt{\alpha(K,x)^2-1}\right)$. However, in
the Bernstein problem the values of $V_K$ are much more of interest for {\em complex} points
$z=x+iy$, in particular for $x\in K$ and $y$ small and nonzero. More precisely, the important
quantity is the normal (sub)derivative
\begin{equation}\label{normalder}
D_{y}^{+} V_E(x):=\lim\inf_{\epsilon\to 0} \frac{V_E(x+i\epsilon y)}{\epsilon}~.
\end{equation}
\begin{oldresult}[\bf Baran \cite{BARAN2}\rm]\label{th:BaranBernstein}
Let $E\subset X$ be any bounded, closed set, $x\in {\rm int} E$ and $0\ne y\in X$. Then for all
$p\in \PPn(X)$ we have
\begin{equation}\label{Baranestimate}
|\langle Dp(x), y \rangle|\le n D_{y}^{+} V_E(x) \sqrt{||p||^2_E-p(x)^2}~.
\end{equation}
\end{oldresult}
Also, one can use the inscribed ellipse method for the estimation of $D_y p(x)$. In the special
case of the standard simplex the yield of both methods can be calculated explicitly. So one can
compare.
\begin{corollary}[\bf Milev, R\'ev\'esz \cite{MR}\rm]\label{cor:coincidence} The  estimate \eqref{Baranestimate}, calculated for the standard simplex $\Delta$ of $\RR^d$ at any point $x\in \Delta$ and in any direction $y\in S^*$ gives exactly identical result to the yield of the inscribed ellipse method.
\end{corollary}
Much remains to explain in this striking coincidence.

There are further yields of the theory of PSH functions, when applied to the Bernstein problem.
For more precise notation now we introduce (interpreting 0/0 as 0 here)
\begin{definition}\label{gradvectorset}
\begin{equation}\label{dradset}
G(E,x):=\{\frac{\nabla p(x)}{n \sqrt{\|p\|^2-p(x)^2}} ~:~ {\bf 0} \neq p \in \PPn, n\in \NN \},
\end{equation}
and following Baran we consider also the convex hull
\begin{equation}\label{conG}
\widetilde{G}(E,x):=\con G(E,x)~.
\end{equation}
\end{definition}
Clearly for any compact $E\subset \RR^d$ $\sup_{n\in\NN} B_n(E,x) = \sup_{u\in G(E,x)} \|u\|$
holds.
\begin{oldresult}[\bf Baran, \cite{BARAN3}\rm]\label{th:equivabove} Let $E$ be a
compact subset of $\RR^d$ with nonempty interior. Then the equilibrium measure $\lambda|_E$ is
absolutely continuous in the interior of $E$ with respect to the Lebesgue measure of $\RR^d$.
Denote its density function by $\lambda(x)$ for all $x\in \intt E$. Then we have $\frac1 {d!}
\lambda(x)\ge \vol \widetilde{G}(E,x)$ for a.a. $x\in\intt E$. Moreover, if $E$ is a symmetric
convex domain of $\RR^d$, then here we have exact equality.
\end{oldresult}

\begin{conjecture}[\bf Baran, \cite{BARAN3}\rm]\label{conj:Baran}
Even if $E$ is a non-symmetric convex body in $\RR^d$ we have $\frac 1{d !} \lambda(x) = \vol
\widetilde{G}(E,x)$.
\end{conjecture}

However, in our recent analysis \cite{R} we found that for dimension
2 $\widetilde{G}(\Delta,x)\subset E_x$ with some ellipsoid domain
$E_x$ of area $\lambda(x)/2$ and major axis \emph{larger} than
\eqref{alphaconjvalue}. So we close this paper with the following
corollary.

\begin{corollary}\label{cor:contradiction} The two conjectures Conjecture \ref{alphaconj} and Conjecture \ref{conj:Baran} can not hold true simultaneously.
\end{corollary}

\section{Acknowledgements}
The above research was either executed or at least started during the authors' stay at the
National Technical University of Athens, Greece in 2001. The author thanks for the support of
the European Commission under the Marie Curie Fellowship contract HPMF-CT-2000-00670.

\end{document}